\documentclass[10pt]{article}
\pdfoutput=1 

\usepackage[utf8]{inputenc}
\usepackage[T1]{fontenc}
\usepackage[english]{babel}

\usepackage{amsmath, amssymb}

\usepackage{geometry}
\usepackage{authblk}

\usepackage{amsmath, amsfonts, amssymb}
\usepackage{bm} 

\usepackage[algo2e, ruled, lined, longend]{algorithm2e}

\usepackage[dvipsnames]{xcolor} 
\usepackage{listofitems}

\usepackage{hyperref}
\usepackage{csquotes}

\usepackage[backend=biber,
            style=numeric-comp,
            sorting=nty,
            giveninits=true,
            maxbibnames=99,
            date=year,      
        hyperref]{biblatex}

\AtEveryBibitem{
  \clearfield{urlyear}
  \clearfield{urlmonth}
  \clearfield{month}
  \clearfield{day}
  \clearfield{eprintclass}
  \clearlist{language}
}
\hypersetup{
    hidelinks,
    colorlinks
}

\usepackage{subcaption}
\usepackage{siunitx}
\usepackage{booktabs}

\usepackage{tikz}
\usetikzlibrary{shapes.geometric,shapes.arrows}
\usetikzlibrary{positioning,fit}

\tikzset{
  nn node/.style={thick,circle,draw=NavyBlue,minimum size=20,inner sep=0.5,outer sep=0.6},
  nn node in/.style={nn node,Green!20!black,draw=Green,fill=Green!25},
  nn node hidden/.style={nn node,NavyBlue!20!black,draw=NavyBlue,fill=NavyBlue!20},
  nn node convol/.style={nn node,orange!20!black,draw=orange,fill=orange!20},
  nn node out/.style={nn node,Red!20!black,draw=Red,fill=Red!20},
  nn connect/.style={thick,Gray}, 
  nn connect arrow/.style={-{Latex[length=4,width=3.5]},thick,NavyBlue,shorten <=0.5,shorten >=1},
  nn node 1/.style={nn node in}, 
  nn node 2/.style={nn node hidden},
  nn node 3/.style={nn node out},
  roundedbox/.style={rectangle,draw,thick,fill=white,outer sep=0pt,align=center,minimum height=3.5em,rounded corners}
}

\tikzset{
  unet txt/.style={text=Gray!60!black, font=\large},
  unet img/.style={rectangle, draw, semithick, inner sep=0pt, outer sep=0pt},
  unet arrow/.style={unet txt, font=\Large, single arrow, draw, semithick, inner sep=3pt, outer sep=0pt, minimum height=40pt},
  unet in/.style={draw=Green, fill=Green!25},
  unet out/.style={draw=Red, fill=Red!25},
  unet enc/.style={draw=NavyBlue, fill=NavyBlue!20},
  unet trans/.style={draw=Cyan, fill=Cyan!20},
  unet updown/.style={draw=orange, fill=orange!20},
  unet dec/.style={draw=Plum, fill=Plum!20},
  unet img in/.style={unet img, unet in},
  unet img out/.style={unet img, unet out},
  unet img enc/.style={unet img, unet enc},
  unet img dec/.style={unet img, unet dec},
  unet img concat/.style={unet img, draw=none, inner sep=0pt},
  unet arrow in/.style={unet arrow, unet trans},
  unet arrow out/.style={unet arrow, unet out},
  unet arrow enc/.style={unet arrow, unet trans, thick, inner ysep=10pt, minimum width=50pt},
  unet arrow updown/.style={unet arrow, unet updown, minimum height=20pt},
  unet arrow bot/.style={unet arrow in},
  unet arrow concat/.style={unet arrow enc},
  unet skip/.style={->, draw=Gray!60!black, line width=3pt},
}

\newlength{\clength}
\newlength{\hzero}
\newlength{\hone}
\newlength{\htwo}
\newlength{\hthree}
\newlength{\hfour}

\newcommand{\imglab}[3]{%
  \node[unet txt, above=0pt of #1] {#2};
  \node[unet txt, left=0pt of #1.south west, anchor=south west, rotate=90, inner xsep=0pt] {#3};
}%

\DeclareMathOperator*{\argmin}{arg\,min}
\DeclareMathOperator{\Span}{span} 
\DeclareMathOperator{\nnz}{nnz} 
\DeclareMathOperator{\conv}{Conv}
\DeclareMathOperator{\convt}{ConvT}
\DeclareMathOperator{\bn}{BN}
\DeclareMathOperator{\prelu}{PreLU}

\newcommand{\NN}{\mathcal{N}_\theta}

\newcommand{\loss}{\mathcal{L}_\theta}

\newcommand{\cost}{J}
\newcommand{\hcost}{\hat{\cost}}

\newcommand{\simiid}{\stackrel{\textnormal{iid}}{\sim}}
\newcommand{\bydef}{:=}

\providecommand{\vec}{}
\renewcommand{\vec}[1]{\mathbf{#1}}

\newcommand{\R}{\mathbb{R}}
\newcommand{\C}{\mathbb{C}}

\newcommand*{\indic}{\text{\usefont{U}{bbold}{m}{n}1}}

\newcommand{\vx}{\vec{x}}

\providecommand{\sp}{}
\renewcommand{\sp}[1]{\tilde{#1}}

\newcommand{\train}{\textnormal{train}}
\newcommand{\val}{\textnormal{val}}
\newcommand{\inner}{\textnormal{inner}}

\newcommand{\pml}{\textnormal{pml}}
\newcommand{\frob}{\textnormal{F}}

\newcommand{\Opml}{\Omega_{\pml}}
\newcommand{\Oinner}{\Omega_{\inner}}

\providecommand{\L}{}
\renewcommand{\L}{L}
\newcommand{\Lpml}{\ell}
\newcommand{\Linner}{\L_{\inner}}

\newcommand{\cin}{c_{\textnormal{in}}}
\newcommand{\cout}{c_{\textnormal{out}}}

\begin{document}

\title{Neural network preconditioning: a case study for the solution of the parametric Helmholtz equation}

\author[1,2]{Luc Giraud}
\author[2]{Carola Kruse}
\author[2,*]{Paul Mycek}
\author[1,3]{Maksym Shpakovych}
\author[1,2,4]{Yanfei Xiang}

\affil[1]{Inria, Concace joint team between Airbus CR\&T, Cerfacs and Inria.}
\affil[2]{Cerfacs, Concace joint team between Airbus CR\&T, Cerfacs and Inria.}
\affil[3]{Current affiliation: Opticalp, France.}
\affil[4]{Current affiliation: IRMA, Université de Strasbourg, MACARON, Inria, France.}
\affil[*]{Corresponding author (\url{mycek@cerfacs.fr}).}

\date{}

\maketitle

{\let\thefootnote\relax\footnotetext{%
Distributed under a \href{https://creativecommons.org/licenses/by/4.0/}{Creative Commons Attribution 4.0} International License.}}

\begin{abstract}
This work presents a hybrid numerical approach for solving linear systems arising from the discretization of the two-dimensional parametric Helmholtz equation. 
A convolutional neural network based on the U-Net architecture is trained in an self-supervised manner (i.e., without explicitly generating reference solutions) to approximate the inverse of the discretized Helmholtz operator, using a loss function involving the residual norm of the linear system. 
The trained network is used as a nonlinear preconditioner within the flexible GMRES (FGMRES) algorithm. 
Numerical experiments show that while the neural network is not accurate enough to act as a standalone solver, it significantly improves the convergence of FGMRES when employed as a preconditioner. 
The neural preconditioner demonstrates robust performance and generalization capabilities with respect to variations in the velocity field and the domain size.
Comparisons with algebraic preconditioners, based on LU factorizations of approximations of the original system matrices, indicate the superior efficiency of the neural approach under equivalent conditions. 
We believe that the proposed method is not tied to a specific neural architecture and can be extended to other parametric PDEs.
\end{abstract}

\section{Introduction}

In this work, we explore a scientific machine learning strategy for solving multiple linear systems arising from the discretization of a parametric PDE.
The approach is based on the hybridization of machine learning techniques with classical numerical algebra methods.
Specifically, a trained neural network (NN) is used as a nonlinear preconditioner for the flexible GMRES (FGMRES) subspace iteration method.
We focus on solving the Helmholtz equation that is used for modeling the propagation of acoustic waves in heterogeneous media.
Specifically, we are interested in the solution of the linear systems resulting from the discretization of this equation in a 2D domain using subspace iterative methods.
For these types of linear systems, no generic, efficient preconditioner currently exists.
In this respect, the objective of this work is to investigate the possibility of using machine learning to design such a preconditioner.
In more classical approaches such as physics-informed neural networks (PINNs)~\cite{karniadakis2021physics}, the trained neural network aims at directly finding the solutions.
The main advantage of our hybrid approach is that the subspace method ensures that the quality of the computed solution far surpasses the quality of the solution predicted by a neural network alone.

Given that neural networks have a nonlinear response to their inputs and are typically trained in 32-bit arithmetic, we consider the FGMRES~\cite{saad:93} subspace method to enable 
a 64-bit computation suited for ill-conditioned problems.
To our knowledge, this idea was first proposed in \cite[Chapter~5]{Xiang2022_SolutionLargeLinear} and further exploited in~\cite{Shpakovych2023_NeuralNetworkPreconditioning}.
We propose to train an NN in a self-supervised manner, using a loss function based on the scaled residual norm of the parametric linear system, for a certain class of parameters.
As a result, the trained NN can be interpreted as a nonlinear approximation of the (linear) inverse of the system matrix, for any parameter instance within the considered parameter class.
We show that the trained NN can thus serve as an effective nonlinear preconditioner for FGMRES with attractive generalization properties.
In particular, the learned NN preconditioner exhibits remarkable robustness when applied to larger problems (in terms of domain size and number of degrees of freedom), or to problems with parameter values or structure outside the training scope.
Our conclusions are not restricted to this network choice and our claim is that one may use their preferred network~\cite{Shpakovych2023_NeuralNetworkPreconditioning} and apply similar ideas to other linear systems.

In recent years, the use of machine learning techniques, particularly NNs, has been considered in various ways for developing preconditioners for Krylov subspace methods.
One contribution in this field, close to ours, comes from Oseledets and colleagues~\cite{rudikov2024neural}, who introduced a discretization-invariant neural operator used as a nonlinear preconditioner within the flexible conjugate gradient (FCG) method, enabling efficient solutions of elliptic PDEs across multiple linear systems.
Dimola et al.~\cite{Dimola2025_AcceleratingSolversMixed,Dimola2025_NeuralPreconditioningKrylov,Dimola2026_NumericalSolutionMixeddimensional} and Nieto Juscafresa~\cite{nieto2024graph} also proposed data-driven preconditioners for (F)GMRES solvers.
In their work, Trifonov et al.~\cite{trifonov2025efficient} proposed a GNN-based preconditioner that learns optimal sparsity patterns and improves the convergence of iterative solvers for complex physical models.
In a related study, they further refined this concept for conjugate gradient solvers, demonstrating significant reduction in iteration counts and computational cost~\cite{trifonov2024learning}. 
Furthermore, Kopaničáková and Karniadakis~\cite{Kopanikov2025} proposed a hybrid preconditioning strategy that leverages DeepONet models to capture low-frequency error components while classical iterative methods handle high-frequency ones. Their framework demonstrated robust performance and convergence acceleration across a wide range of parameterized linear systems, offering a novel perspective distinct from purely GNN-based or multigrid techniques. This approach complements other neural preconditioning methods and reveals the potential of operator learning frameworks in the preconditioning landscape.
Yusuf et al.~\cite{yusuf2024constructing}, who designed incomplete LU (ILU) preconditioners via GNNs for advection-dominated PDEs; and Li et al.~\cite{li2025neural}, who introduced neural operators to accelerate PDE solvers directly.
Azulay and Treister~\cite{Azulay2023_MultigridAugmentedDeepLearning} combined multigrid methods with deep learning for Helmholtz problems, while Moore~\cite{moore2022machine} explored residual neural networks for constructing adaptive preconditioners.

The field of scientific machine learning is rapidly evolving, and the few references provided herein are not intended to offer an exhaustive overview of the state of the art, as such an undertaking would far exceed the scope of the present work.

\section{Background}

\subsection{The Helmholtz equation}\label{sec:helm}
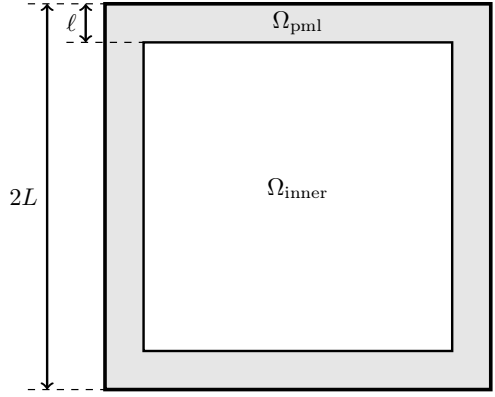
\begin{figure}[!ht]\centering%
\resizebox{.4\textwidth}{!}{%
\begin{tikzpicture}[x=8pt,y=8pt]%
\node[rectangle, draw=black, line width=1.5pt, fill=gray!20, fit={(0,0)(20,20)}, inner sep=0pt] {};
\node[rectangle, draw=black, line width=1pt, fill=white, fit={(2,2)(18,18)}, inner sep=0pt] {$\Oinner$};
\node at (10, 19) {$\Opml$};
\draw[dashed, line width=.5pt] (-4,20) -- (2,20);
\draw[dashed, line width=.5pt] (-2,18) -- (2,18);
\draw[dashed, line width=.5pt] (-4,0) -- (2,0);
\draw[<->, line width=1pt] (-1,18) -- node[midway, left]{$\ell$} (-1,20) ;
\draw[<->, line width=1pt] (-3,0) -- node[midway, left]{$2L$} (-3,20); 
\end{tikzpicture}%
}%
\caption{Illustration of the PML.}
\label{fig:pml}
\end{figure}

We consider the two-dimensional Helmholtz equation, subject to the so-called Sommerfeld radiation condition,
\begin{equation}\label{eq:helmholtz}
    \begin{cases}
        \nabla^2 u + k^2 u = f \; \mbox{ in } \R^2,\\[\medskipamount]
        \lim_{\|\vx\|_2 \to \infty} \|\vx\|_2^{1/2} 
        \left(
            \dfrac{\partial u}{\partial \|\vx\|_2}(\vx) - \jmath 
            k(\vx)
            u(\vx)
        \right)
        = 0,
    \end{cases}
\end{equation}
where $\jmath \in \C$ denotes the imaginary unit such that $\jmath^2=-1$, and $u \colon \R^2 \to \C$, $f \colon \R^2 \to \C$, $k \colon \R^2 \to \R^+$ are scalar fields, referred to as the solution field, the source field, and the wavenumber field, respectively.
The latter is defined by $k \colon \vx \mapsto \omega / c(\vx)$, where $\omega \in \R^+$ is the angular frequency of the source, and $c \colon \R^2 \to \R^+$ is the velocity field.
In what follows, we restrict ourselves to the  case $\omega = 1$.
Numerically, the unbounded domain $\R^2$ is truncated to a computational, square domain $\Omega \bydef [-\L, \L]^2 \subset \R^2$ and
the derivatives in \eqref{eq:helmholtz} are discretized using Fourier differentiation~\cite{Stanziola2021_HelmholtzEquationSolver, Xiang2022_SolutionLargeLinear, Shpakovych2023_NeuralNetworkPreconditioning} with $N$ points in each direction.
Consequently, there is a total of $n \bydef N^2$ degrees of freedom.
The boundary conditions are (approximately) enforced using perfectly matched layers (PMLs)~\cite{Berenger1994_PerfectlyMatchedLayer} in the outer region $\Opml \bydef \Omega \setminus \Oinner$, where $\Oinner \bydef [-\Linner, \Linner]^2$, with $\Linner \bydef \L - \Lpml$ and $\Lpml < \L$, as illustrated in Figure~\ref{fig:pml}.

The discretization of the Helmholtz problem~\eqref{eq:helmholtz} yields a system of linear equations of the form $A(c)x = b$, involving the discrete velocity field $c \in \R^{n}$ and the associated discretized operator (matrix) $A(c) \in \C^{n\times n}$, as well as the discretized source term $b \in \C^{n}$ and solution $x \in \C^{n}$.

\subsection{Subspace solver with nonlinear preconditioner}\label{sec:nla}

The solution of large linear systems of equations is often approached by iterative methods.
Among the advantages of these techniques are their moderate memory consumption and the ability to stop the iterations when the quality of the solution is similar to the possible uncertainty in the matrix entries or right-hand sides, which may arise from discretization errors in the case of a PDE solution or from uncertainty in some input data. Backward error analysis~\cite{wilkinson1963rounding, Rigal1967} provides powerful techniques for designing meaningful stopping criteria.
In recent decades, subspace methods, such as Krylov subspace methods, have attracted much attention and have become standard techniques for solving linear systems. For unsymmetric problems, the best known method is certainly GMRES~\cite{sasc:86}, which computes the minimum residual norm approximation in a Krylov subspace. This technique relies on the so-called Arnoldi process, which incrementally builds an orthonormal basis of the nested subspaces. At step $k$, the Arnoldi algorithm computes a set of orthonormal vectors $V_k = [v_1, \dots, v_k ]$ that satisfy the so-called Arnoldi relation, which is written in matrix form:%
\begin{equation}%
A V_k = V_{k+1}\bar{H}_k \quad \mbox{\rm with } V_{k+1}^H V_{k+1} = I_{k+1}
\end{equation}%
where $\bar{H}_k \in \mathbb{C}^{(k+1)\times k}$ is upper Hessenberg.
The minimum residual norm iterate $x_k \in \Span\{v_1, \dots, v_k \}$, that is $x_k = V_k y_k$ with  $y_k =\argmin_{ y \in \C^k} (\| \bar{H}_k y - \beta e_1 \|_2 )$ where $\beta =\| b \|_2$.
In practice, a preconditioner $M$ is used to speed up the convergence of GMRES, where $M$ is expected to approximate $A^{-1}$ somehow. For GMRES it is recommended to use the preconditioner on the right, i.e., GMRES solves $AM t = b$, with $Mt=x$. Thus, at each iteration, GMRES still minimizes the residual norm of the residual associated with the original linear system. The Arnoldi relation is%
\begin{equation}
A M V_k = V_{k+1}\bar{H}_k \quad \mbox{\rm with } V_{k+1}^H V_{k+1} = I_{k+1}.
\end{equation}
In~\cite{saad:93}, Saad introduced the idea of a flexible preconditioner, which allows to have a different preconditioning matrix $M_i$ at each iteration, so that the generalized Arnoldi relation becomes
\begin{equation}\label{eq:generalizedArnoldi}
    A Z_k = V_{k+1}\bar{H}_k \quad \mbox{\rm with } V_{k+1}^H V_{k+1} = I_{k+1}, 
\end{equation}
where
$Z_k = [z_1, \dots, z_k]$ with $z_i = M_i v_i~(i=1,\dots,k)$.
This idea can easily be extended to the situation where the preconditioner is no longer linear, i.e., $z_i = {\cal M}_i(v_i)$, where each ${\cal M}_i$ may be a nonlinear operator. 
In the latter case, Equation~\eqref{eq:generalizedArnoldi} still holds and uniquely defines each iterate as long as $Z_k$ remains full rank.
In this paper we use this equality with ${\cal M}_i() = \NN()$, which is the trained neural network (see Section~\ref{sec:nn-components}). The nonlinearity has two origins; first, the activation functions in the neurons are nonlinear, and second, the preconditioner is computed in 32-bit arithmetic while the rest of the computation is performed in 64-bit arithmetic. This truncation due to casting 64 bits to 32 bits is also nonlinear.
A sketch of the NN-preconditioned FGMRES algorithm is provided in Algorithm~\ref{alg:nn-fgmres}, where the key step of applying the NN preconditioner $\NN()$ is performed on line~\ref{algline:fgmres-nn-prec}.

{\DontPrintSemicolon\LinesNumbered%
\begin{algorithm2e}[ht]%
  $r_0 = b - A x_0$, $\beta = \|r_0\|_2$, $v_1 = r_0 / \beta$ \;
  \For{$j=1,\ldots,m$}{
    $z_j = \NN(v_j)$ \label{algline:fgmres-nn-prec} \;
    $w = A z_j$ \; 
    \For{$i=1,\ldots,j$}{
      $h_{i,j} = \langle w , v_i \rangle$\;
      $w \leftarrow w - h_{i,j} v_i$ \;
    }
    $h_{j+1, j} = \|w\|_2$ \;
    $v_{j+1} = w / h_{j+1, j}$ \;
  }
  $\bar{H}_m = [h_{i,j}]_{1\leq j \leq m, 1 \leq i \leq j+1}$ \;
  $Z_m = [z_1, \dots, z_m]$ \;
  $y_m = \argmin_{y \in \C^m} \| \bar{H}_m y - \beta e_1 \|_2$ \;
  $x_m = x_0 + Z_m y_m$ \;
  \caption{Sketch of the NN-preconditioned FGMRES algorithm.}
  \label{alg:nn-fgmres}
\end{algorithm2e}%
}

\section{Neural network components}\label{sec:nn-components}

\subsection{Network architecture}\label{sec:network-architecture}

As stated in Section~\ref{sec:nla}, we propose to use a neural network as a nonlinear preconditioner.
To do so, we design a neural network based on the U-Net architecture~\cite{Ronneberger2015_UNetConvolutionalNetworks}, which takes the discretized source and velocity fields as inputs.
To make them suitable to the convolutional nature of the U-Net, these 2D fields are recast as image-like tensors.
Specifically, recalling that $n=N^2$, the complex-valued source term $b \in \C^n$ may be recast as a $2 \times N \times N$ tensor composed of two $N \times N$ channels representing the real and imaginary parts of the field.
Similarly, the real-valued velocity field $c \in \R^n$ may be recast as a $1 \times N \times N$ tensor.
These two inputs are combined in a single $3 \times N \times N$ input tensor.
The output of the network is a $2 \times N \times N$ tensor, representing a 2D, complex-valued field.
As we shall see in Section~\ref{sec:training}, the network will be trained to provide approximations (predictions) $\hat{x}_\theta$ of discrete solution fields $x$ of the linear system $A(c)x = b$ corresponding to a discretized version of the Helmholtz problem \eqref{eq:helmholtz}.

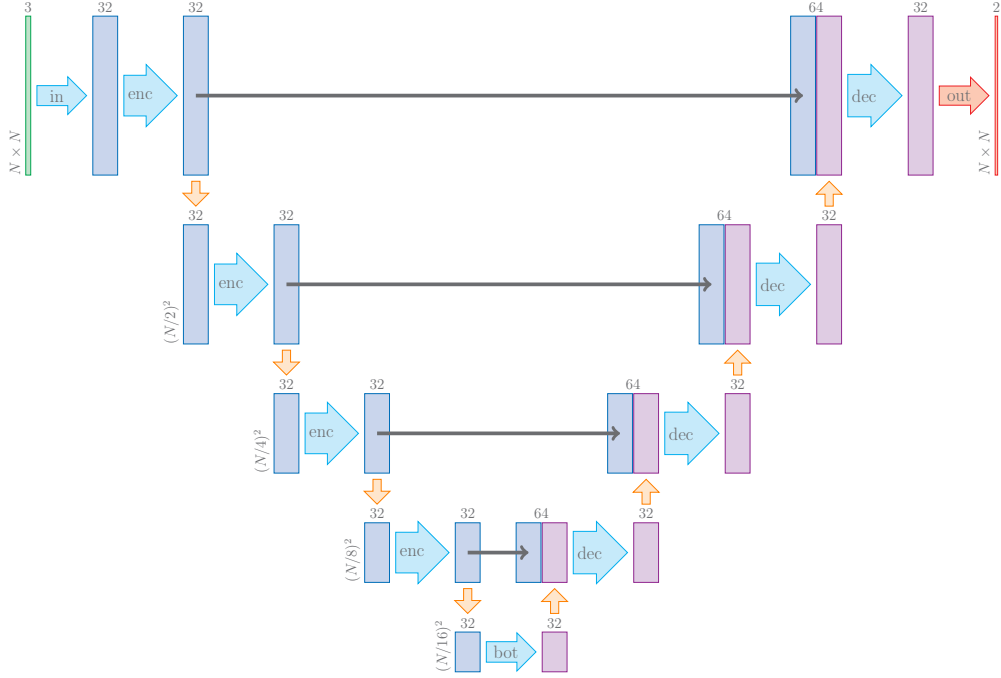
\begin{figure}[!ht]\centering%
\resizebox{0.8 \textwidth}{!}{%
\begin{tikzpicture}%
\node[unet img in, minimum width=4pt, minimum height=\hzero] (img-in) {};
\imglab{img-in}{3}{$N \times N$}
\node[unet arrow in, right=5pt of img-in] (arr-in) {in};
\node[unet img enc, minimum width=\clength, minimum height=\hzero, right=5pt of arr-in] (img-enc-00) {};
\imglab{img-enc-00}{32}{}
\node[unet arrow enc, right=5pt of img-enc-00] (arr-trans-0) {enc};
\node[unet img enc, minimum width=\clength, minimum height=\hzero, right=5pt of arr-trans-0] (img-enc-01) {};
\imglab{img-enc-01}{32}{}
\node[unet arrow updown, below=5pt of img-enc-01.south, anchor=west, rotate=-90] (arr-down-0) {};
\node[unet img enc, minimum width=\clength, minimum height=\hone, below=15pt of arr-down-0.east, anchor=north] (img-enc-10) {};
\imglab{img-enc-10}{32}{$(N/2)^2$}
\node[unet arrow enc, right=5pt of img-enc-10] (arr-trans-1) {enc};
\node[unet img enc, minimum width=\clength, minimum height=\hone, right=5pt of arr-trans-1] (img-enc-11) {};
\imglab{img-enc-11}{32}{}
\node[unet arrow updown, below=5pt of img-enc-11.south, anchor=west, rotate=-90] (arr-down-1) {};
\node[unet img enc, minimum width=\clength, minimum height=\htwo, below=15pt of arr-down-1.east, anchor=north] (img-enc-20) {};
\imglab{img-enc-20}{32}{$(N/4)^2$}
\node[unet arrow enc, right=5pt of img-enc-20] (arr-trans-2) {enc};
\node[unet img enc, minimum width=\clength, minimum height=\htwo, right=5pt of arr-trans-2] (img-enc-21) {};
\imglab{img-enc-21}{32}{}
\node[unet arrow updown, below=5pt of img-enc-21.south, anchor=west, rotate=-90] (arr-down-2) {};
\node[unet img enc, minimum width=\clength, minimum height=\hthree, below=15pt of arr-down-2.east, anchor=north] (img-enc-30) {};
\imglab{img-enc-30}{32}{$(N/8)^2$}
\node[unet arrow enc, right=5pt of img-enc-30] (arr-trans-3) {enc};
\node[unet img enc, minimum width=\clength, minimum height=\hthree, right=5pt of arr-trans-3] (img-enc-31) {};
\imglab{img-enc-31}{32}{}
\node[unet arrow updown, below=5pt of img-enc-31.south, anchor=west, rotate=-90] (arr-down-3) {};
\node[unet img enc, minimum width=\clength, minimum height=\hfour, below=15pt of arr-down-3.east, anchor=north] (img-enc-40) {};
\imglab{img-enc-40}{32}{$(N/16)^2$}
\node[unet arrow bot, right=5pt of img-enc-40] (arr-bot) {bot};
\node[unet img dec, minimum width=\clength, minimum height=\hfour, right=5pt of arr-bot] (img-dec-41) {};
\imglab{img-dec-41}{32}{}
\node[unet arrow updown, above=15pt of img-dec-41.north, anchor=west, rotate=90] (arr-up-3) {};
\node[unet img dec, minimum width=\clength, minimum height=\hthree, above=5pt of arr-up-3.east, anchor=south] (img-dec-30) {};
\node[unet img enc, minimum width=\clength, minimum height=\hthree, left=1pt of img-dec-30] (img-skip-3) {};
\node[unet img concat, fit=(img-dec-30) (img-skip-3)] (img-concat-3) {};
\draw[unet skip] (img-enc-31.center) -- (img-skip-3.center);
\node[unet txt, above=0pt of img-concat-3] {64};
\node[unet arrow concat, right=5pt of img-concat-3] (arr-trans-dec-3) {dec};
\node[unet img dec, minimum width=\clength, minimum height=\hthree, right=5pt of arr-trans-dec-3] (img-dec-31) {};
\imglab{img-dec-31}{32}{}
\node[unet arrow updown, above=15pt of img-dec-31.north, anchor=west, rotate=90] (arr-up-2) {};
\node[unet img dec, minimum width=\clength, minimum height=\htwo, above=5pt of arr-up-2.east, anchor=south] (img-dec-20) {};
\node[unet img enc, minimum width=\clength, minimum height=\htwo, left=1pt of img-dec-20] (img-skip-2) {};
\node[unet img concat, fit=(img-dec-20) (img-skip-2)] (img-concat-2) {};
\draw[unet skip] (img-enc-21.center) -- (img-skip-2.center);
\node[unet txt, above=0pt of img-concat-2] {64};
\node[unet arrow concat, right=5pt of img-concat-2] (arr-trans-dec-2) {dec};
\node[unet img dec, minimum width=\clength, minimum height=\htwo, right=5pt of arr-trans-dec-2] (img-dec-21) {};
\imglab{img-dec-21}{32}{}
\node[unet arrow updown, above=15pt of img-dec-21.north, anchor=west, rotate=90] (arr-up-1) {};
\node[unet img dec, minimum width=\clength, minimum height=\hone, above=5pt of arr-up-1.east, anchor=south] (img-dec-10) {};
\node[unet img enc, minimum width=\clength, minimum height=\hone, left=1pt of img-dec-10] (img-skip-1) {};
\node[unet img concat, fit=(img-dec-10) (img-skip-1)] (img-concat-1) {};
\draw[unet skip] (img-enc-11.center) -- (img-skip-1.center);
\node[unet txt, above=0pt of img-concat-1] {64};
\node[unet arrow concat, right=5pt of img-concat-1] (arr-trans-dec-1) {dec};
\node[unet img dec, minimum width=\clength, minimum height=\hone, right=5pt of arr-trans-dec-1] (img-dec-11) {};
\imglab{img-dec-11}{32}{}
\node[unet arrow updown, above=15pt of img-dec-11.north, anchor=west, rotate=90] (arr-up-0) {};
\node[unet img dec, minimum width=\clength, minimum height=\hzero, above=5pt of arr-up-0.east, anchor=south] (img-dec-00) {};
\node[unet img enc, minimum width=\clength, minimum height=\hzero, left=1pt of img-dec-00] (img-skip-0) {};
\node[unet img concat, fit=(img-dec-00) (img-skip-0)] (img-concat-0) {};
\draw[unet skip] (img-enc-01.center) -- (img-skip-0.center);
\node[unet txt, above=0pt of img-concat-0] {64};
\node[unet arrow concat, right=5pt of img-concat-0] (arr-trans-dec-0) {dec};
\node[unet img dec, minimum width=\clength, minimum height=\hzero, right=5pt of arr-trans-dec-0] (img-dec-01) {};
\imglab{img-dec-01}{32}{}
\node[unet arrow out, right=5pt of img-dec-01] (arr-out) {out};
\node[unet img out, minimum width=2pt, minimum height=\hzero, right=5pt of arr-out] (img-out) {};
\imglab{img-out}{2}{$N \times N$}
\end{tikzpicture}%
%
}%
\caption{Overview of our U-Net architecture, whose main building blocks are described in Tables~\ref{tab:unet_seq_operations} and~\ref{tab:unet_unit_operations}.}
\label{fig:Unet}
\end{figure}%

The specific U-Net architecture prescribed in this work is depicted in Figure~\ref{fig:Unet}.
The modules represented in Figure~\ref{fig:Unet} are described in Table~\ref{tab:unet_seq_operations}.
They consist of sequences of parameterized layers detailed in Table~\ref{tab:unet_unit_operations}.
The downscaling convolutions, represented by the downward orange arrows, are designed such that the size of the ``image'' is halved.
Similarly, the upscaling convolutions, represented by the upward orange arrows, are designed such that the size of the ``image'' is multiplied by 2.
As reported in Table~\ref{tab:unet_seq_operations}, our U-Net architecture comprises a total of \num{831650} trainable parameters, which is independent of $n$ owing to its convolutional nature.

\begin{table}[!ht]
    \centering\small
    \begin{tabular}{>{\centering}m{1cm}lccc}
    \toprule
    \bfseries Module & \bfseries Sequence of layers & \bfseries \#params/module & \bfseries \#modules & \bfseries total \#params\\
    \midrule
    \resizebox{\linewidth}{!}{\tikz{\node[unet arrow in] {in};}} &  
        $\conv_{3, 32, 3, 1, 1}, \bn_{32}, \prelu_{32}$ & \num{992} & 1 & \num{992}\\[1em]
    \resizebox{\linewidth}{!}{\tikz{\node[unet arrow in] {bot};}} &  
        $\conv_{32, 32, 3, 1, 1}, \bn_{32}, \prelu_{32}$ & \num{9344} & 1 & \num{9344}\\[1em]
    \resizebox{\linewidth}{!}{\tikz{\node[unet arrow enc] {enc};}} &
        $\conv_{32, 32, 3, 1, 1}, \bn_{32}, \prelu_{32},\conv_{32, 32, 3, 1, 1}$ & \num{18592} & 4 & \num{74368} \\[2em]
    \resizebox{\linewidth}{!}{\tikz{\node[unet arrow enc] {dec};}} &
        $\conv_{64, 64, 3, 1, 1}, \bn_{64}, \prelu_{64},\conv_{64, 32, 3, 1, 1}$ & \num{55584} & 4 & \num{222336}\\[2em]
    \resizebox{\linewidth}{!}{\tikz{\node[unet arrow out] {out};}} &
        $\conv_{32, 2, 1, 0, 1}$ & 66 & 1 & 66\\[1em]
    \scalebox{.8}{\tikz{\node[unet arrow updown, rotate=-90] {};}} &  
        $\conv_{32, 32, 8, 3, 2}$ & \num{65568} & 4 & \num{262272} \\[1em]
    \scalebox{.8}{\tikz{\node[unet arrow updown, rotate=90] {};}} &  
        $\convt_{32, 32, 8, 3, 2}$ & \num{65568} & 4 & \num{262272} \\[1em]
    \scalebox{.8}{\tikz{\draw[unet skip] (0,0) -- (1,0);}} & skip connection  & 0 & 4 & 0 \\
    \midrule
    & & \multicolumn{2}{r}{\bfseries Total} & \bfseries\num{831650}\\
    \bottomrule
    \end{tabular}
    \caption{Description of the modules used in the U-Net architecture depicted in Figure~\ref{fig:Unet}.
      Each module is defined as a sequence of layers detailed in Table~\ref{tab:unet_unit_operations}, and the associated number of parameters is deduced from the values provided therein.}
    \label{tab:unet_seq_operations}
\end{table}

\begin{table}[!ht]
    \centering\small
    \begin{tabular}{lllll}
        \toprule    
        \bfseries Layer & \bfseries PyTorch name & \bfseries Arguments & \bfseries \#params & \bfseries Output size  \\
        \midrule
        $\conv_{\cin, \cout, k, p, s}$ & \texttt{Conv2d} & $\texttt{in\_channels} = \cin$ & $(k^2\cin+1)\cout$ & $o=\lfloor (i + 2p - k) / s + 1 \rfloor$ \\
        & & $\texttt{out\_channels} = \cout$ & & \\
        & & $\texttt{kernel\_size} = k$ & & \\
        & & $\texttt{padding} = p$ & & \\
        & & $\texttt{stride} = s$ & & \\
        \midrule
        $\convt_{\cin, \cout, k, p, s}$ & \texttt{ConvTranspose2d} & $\texttt{in\_channels} = \cin$ & $(k^2\cin+1)\cout$ & $o=s (i-1) -2p + k$ \\
        & & $\texttt{out\_channels} = \cout$ & & \\
        & & $\texttt{kernel\_size} = k$ & & \\
        & & $\texttt{padding} = p$ & & \\
        & & $\texttt{stride} = s$ & & \\
        \midrule
        $\bn_{c}$ & \texttt{BatchNorm2d} & $\texttt{num\_features} = c$ & $2c$ & $o=i$ \\
        \midrule
        $\prelu_{c}$ & \texttt{PReLU} & $\texttt{num\_parameters} = c$ & $c$ & $o=i$ \\
        \bottomrule
    \end{tabular}
    \caption{Description and PyTorch name of the NN layers used in the U-Net architecture described in Figure~\ref{fig:Unet} and Table~\ref{tab:unet_seq_operations}, along with their arguments, and the induced number of parameters and output size.}
    \label{tab:unet_unit_operations}
\end{table}

The abstract, functional representation of the neural network with parameters $\theta$ is thus a nonlinear operator,
\begin{equation}
  \label{eq:nn_operator}
  \begin{aligned}
    \NN \colon \C^n \times \R^n & \to \C^n,\\
    (b, c) & \mapsto \NN(b, c) = \NN(c)(b),
  \end{aligned}
\end{equation}
where the conversions to and from image-like tensors are implicit.

\subsection{Training of the network}\label{sec:training}

\begin{figure}[!ht]\centering%
\resizebox{.8\textwidth}{!}{%
\begin{tikzpicture}%
  \node[rectangle, draw, inner ysep=10pt] (nn) {$\NN$};
  \node[left=16pt of nn, yshift=8pt] (b) {$b$};
  \node[left=16pt of nn, yshift=-8pt] (c) {$c$};
  \node[right=16pt of nn] (uhat) {$\hat{x}_\theta = \NN(b,c)$};
  \node[rectangle, draw, inner ysep=6pt, right=16pt of uhat] (A) {$A(\cdot)$};
  \node[right=16pt of A] (bhat) {$\hat{b}_\theta = A(c)\hat{x}_\theta$};
  \node[rectangle, draw, inner ysep=6pt, right=16pt of bhat] (loss) {loss};
  \node[right=16pt of loss] (L) {$\loss=\dfrac{\|b - \hat{b}_\theta \|_2^2}{\|b\|_2^2}$};

  \draw[->] (b) -- (b-|nn.west);
  \draw[->] (c) -- (c-|nn.west);
  \draw[->] (nn) -- (uhat);
  \draw[->] (uhat) -- (A);
  \draw[->] (A) -- (bhat);
  \draw[->] (bhat) -- (loss);
  \draw[->] (loss) -- (L);
  \draw[->] (b) -- ++(0pt,16pt) coordinate (bb) -- (bb-|loss) -- (loss);
  \draw[->] (c) -- ++(0pt,-16pt) coordinate (cc) -- (cc-|A) -- (A);
\end{tikzpicture}%
}%
\caption{Schematic representation of the training framework.}
\label{fig:nn_general_schamatic}
\end{figure}
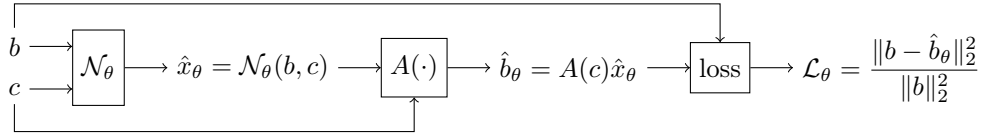

The training process aims at finding the parameters $\theta^* \in \Theta$ of the neural network that minimize a certain cost functional $\cost \colon \Theta \to \R$, where $\Theta$ denotes the parameter space, i.e.,
\begin{equation}
  \label{eq:nn_optim}
  \theta^* \in \argmin_{\theta \in \Theta} \cost(\theta).
\end{equation}
For our specific problem, the cost function takes the form
\begin{equation}
  \label{eq:nn_cost_fn}
  \cost(\theta) = \mathbb{E}_{b,c} [\loss(b,c)],
\end{equation}
where $b$ and $c$ are the inputs of the network, $\loss$ is a loss function and $\mathbb{E}_{b,c}$ denotes the expectation operator with respect to $b$ and $c$.
In this work, we consider the loss function defined by
\begin{equation}
  \label{eq:nn_loss_fn}
  \loss(b,c)
  =
  \dfrac{\|b - A(c)\NN(b,c)\|_2^2}{\|b\|_2^2}
  =
  \dfrac{\|b - A(c)\hat{x}_\theta \|_2^2}{\|b\|_2^2},
\end{equation}
where $\hat{x}_\theta \bydef \NN(b,c) = \NN(c)(b)$ denotes the prediction of the neural network.
The neural network is thus trained to be a good approximation of $A^{-1}$, i.e., roughly speaking, so that $\NN(c) \approx A^{-1}(c)$ for any $c$.
The training process can here be interpreted as being self-supervised, in the sense that the prediction $\hat{x}_\theta$ are not (directly) compared to explicitely-generated reference values (or labels, in machine learning terminology) $x^*$ in the loss function.
We remark that the residual of the linear system for a given prediction $\hat{x}_\theta = \NN(c)(b)$, namely $b - A(c)\hat{x}_\theta$, which is at the core of the loss definition, corresponds to an algebraic (discretized) version of the residual of the PDE, including the boundary conditions.
As such, the loss may be interpreted as being physics-informed.
In terms of numerical linear algebra, we further note that for any given $c \in \R^n$, $\loss(b,c) = \eta_b(\hat{x}_\theta)^2$, where $\eta_b(\tilde{x})$ denotes the backward error with respect to $b$ of $\tilde{x}$ as an approximate solution to the linear system $A(c)x = b$~\cite{higham2002accuracy,Rigal1967}.
An alternative interpretation can be drawn from rewriting the loss function \eqref{eq:nn_loss_fn} as
\begin{equation}
  \label{eq:nn_loss_fn_b}
  \loss(b,c)
  =
  \dfrac{\|b - \hat{b}_\theta \|_2^2}{\|b\|_2^2},
\end{equation}
with $\hat{b}_\theta = A(c)\hat{x}_\theta = A(c)\NN(c)(b)$.
This form of the loss function highlights that the neural network is trained so that $A(c)\NN(c)$ is a good approximation of the identity operator, similar to an autoencoder~\cite{Bishop2024_Autoencoders}.
The neural network $\NN$ can then be interpreted as an encoder, trained such that it is decoded by $A$, although here, unlike conventional autoencoder systems, $\NN$ does not reduce the dimension of the input space.
The training process is summarized in Figure~\ref{fig:nn_general_schamatic}.

\begin{table}[!ht]
  \centering\small
  \begin{tabular}{lc}
    \toprule  
    \textbf{Parameter} & \textbf{Value} \\
    \midrule
    Optimizer & Adam~\cite{Kingma2017_AdamMethodStochastic} \\
    Mini-batch size $B$ & 32 \\
    Number of epochs & 500 \\
    Size of data set $K$ & \num{16000} \\
    Training/validation data set split & 96\%/4\% \\
    Learning rate (LR) & \num{e-3} \\
    Gradient clipping (norm) & 1 \\
    LR scheduler (LRS) & \texttt{ReduceLROnPlateau} \\
    LRS factor & 0.5 \\
    LRS patience & 10 \\
    LRS min LR & \num{e-5} \\
    \bottomrule
  \end{tabular}
  \caption{Parameters of the training process.}
  \label{tab:nn_optim_params}
\end{table}

Numerically, the solution of the optimization problem \eqref{eq:nn_optim} is approximated using the Adam algorithm~\cite{Kingma2017_AdamMethodStochastic} over 500 epochs, from a data set of $K=\num{16000}$ synthetic pairs of inputs $\mathcal{X} = \{(b^{(k)}, c^{(k)})\}_{k=1}^K$, which are cheaply generated from the (prescribed) probability distribution of $(b,c)$ (see below).
The data set is split into a training data set $\mathcal{X}_{\train} \subset \mathcal{X}$ of size \num{15360} (96\% of the data set), on which the training is performed, and a validation data set $\mathcal{X}_{\val} \bydef \mathcal{X} \setminus \mathcal{X}_{\train}$ of size 640 (4\% of the data set), on which the model is validated at the end of each epoch.
The parameters of the optimizer are summarized in Table~\ref{tab:nn_optim_params}.
At each iteration (gradient step) of the Adam algorithm, the cost function is approximated from a mini-batch $\mathcal{B} \subset \mathcal{X}_{\train}$ of size $B=32$ as
\begin{equation}\label{eq:nn_loss_minibatch}
  \cost
  \approx
  \hcost_{\mathcal{B}}
 \bydef 
 \dfrac{1}{B} \sum_{(b,c) \in \mathcal{B}} \loss(b,c).
\end{equation}
The network was trained using training and validation data generated with the following parameters involved in the discretisation of the Helmholtz equation, as introduced in Section~\ref{sec:helm}: $N=64$, $L=20$ and $\ell=4$.
Furthermore, the distribution of $(b,c)$ is defined component-wise and independently for $b$ and $c$ such that the entries of $b$ are independent and normally distributed with zero mean and unit variance (we denote $b \simiid {\cal N}(0,1)$), the ones of $c$ are independent and uniformly distributed between 1 and 2 (i.e., $c \simiid {\cal U}([1,2])$).
We note that the latter implies that the wavelength $\lambda \bydef 2\pi/k = 2\pi c/\omega$ lies within the range $[2\pi, 4\pi]$, so that there are between 10 and 20 discretization points per wavelength, since $2L/N = 5/8$.
These experiments are presented subsequently in this paper.

\section{Numerical experiments}\label{sec:num}

\begin{figure}[!ht]\centering%
    \begin{subfigure}{.49\textwidth}\centering%
    \includegraphics[width=\textwidth]{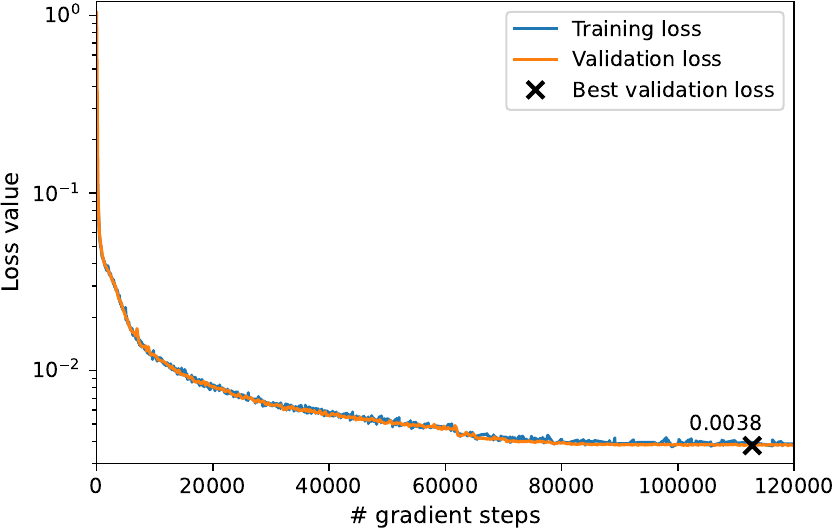}%
    \caption{Training \#1}
    \end{subfigure}\hfill%
    \begin{subfigure}{.49\textwidth}\centering%
    \includegraphics[width=\textwidth]{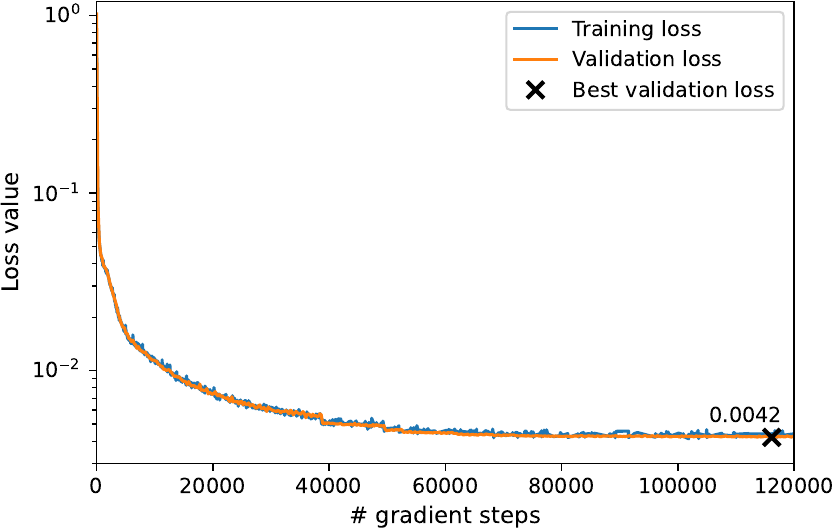}%
    \caption{Training \#2}
    \end{subfigure}%
    \caption{Evolution of the training and validation loss throughout the training process (reported against the number of gradient steps).}
    \label{fig:nn-losses}
\end{figure}

In this section, we perform numerical experiments to evaluate the quality and robustness of the neural network, first as a standalone solver (in other words, as a surrogate), then as a nonlinear preconditioner as described in Section \ref{sec:nla}.
For that purpose we select two trained neural networks that have reached a comparable value of the validation loss function during their training; they are referred to as Training~\#1 and Training~\#2.
The evolution of the training and validation loss throughout each training is reported in Figure~\ref{fig:nn-losses}.
The loss function is defined in Eq.~\eqref{eq:nn_loss_fn} and the training loss value is evaluated as described in Eq.~\eqref{eq:nn_loss_minibatch}; the validation loss value is evaluated similarly from the validation dataset.
We observe that both trainings evolve similarly, that is, the training and validation losses decrease with gradient steps (or, equivalently, with epochs), and the validation loss of both trainings reaches a similar minimum of about $0.004$ after a similar number of gradient steps (or epochs).
We remark that, according to this metric, the network resulting from Training~\#1 appears to be slightly better.

For all the numerical experiments presented in the rest of this paper, the quality of an approximate solution $\tilde{x}$ to the linear system $A(c)x=b$ will be evaluated through the normwise backward error~\cite{higham2002accuracy,Rigal1967} with respect to the right-hand side $b$,
\begin{equation}\label{eq:backward_err_b}
    \eta_b(\tilde{x})
    \bydef
    \dfrac{\|b - A(c)\tilde{x} \|_2}{\|b\|_2}.
\end{equation}
The approximate solution will arise either from the direct use of the NN as a surrogate model serving as a standalone solver (see Section~\ref{sec:standalone}), in which case $\tilde{x} = \NN(b,c)$, or from iterations of the NN-preconditioned FGMRES solver (see Section~\ref{sec:nn-precond}).

\subsection{Numerical behavior of the trained network as a standalone solver}
\label{sec:standalone}

We first investigate the numerical robustness of the two trained networks (corresponding to Training~\#1 and Training~\#2) used as surrogate models that serve as standalone solvers.
We study their quality and robustness both for problems similar to those used during the training, that is, linear systems arising from the discretization on a $64 \times 64$ grid and $c \simiid {\cal U}([1,2])$ as well as for other distributions for $c$ and larger grid sizes.

\subsubsection{Generalization capabilities with respect to the parameter $c$}\label{sec:raw_nn_gen_c}

\begin{figure}[!ht]\centering%
  \begin{subfigure}{.49\textwidth}%
    \includegraphics[width=\textwidth]{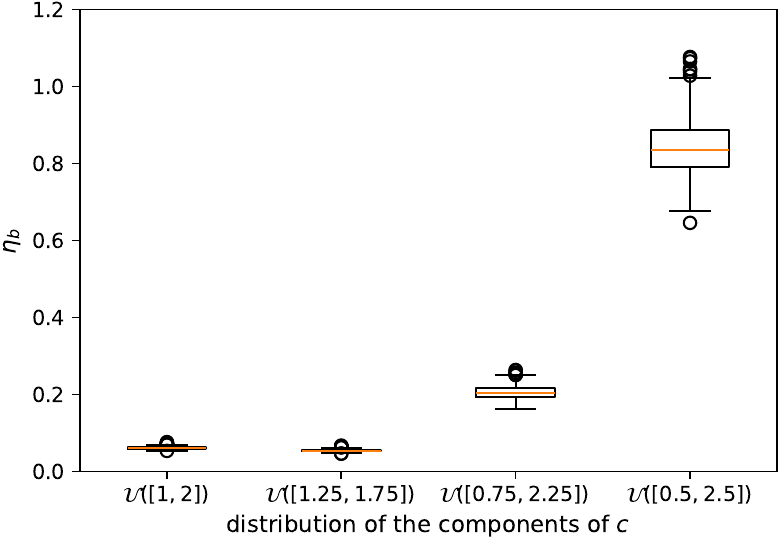}%
    \caption{Training \#1}%
  \end{subfigure}\hfill%
  \begin{subfigure}{.49\textwidth}%
    \includegraphics[width=\textwidth]{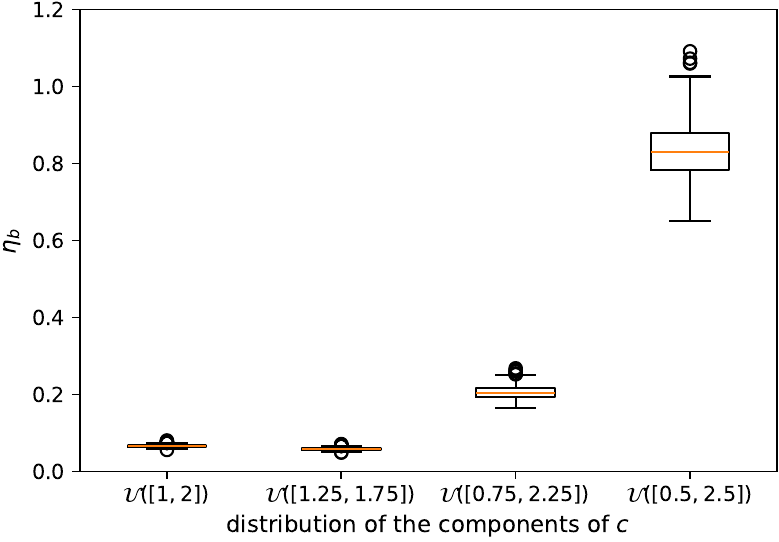}%
    \caption{Training \#2}%
  \end{subfigure}%
  \caption{Direct prediction from the neural network on problems with $c$ having different distributions, where the neural network has been trained with $\alpha=1$ and $c \simiid \mathcal{U}([1,2])$.
    Boxplots of $\eta_b$ based on \num{1000} experiments.}
  \label{fig:raw_nn_gen_sos}
\end{figure}
  
For the set of trained parameters, we display in Figure~\ref{fig:raw_nn_gen_sos} the boxplots of $\eta_b(\hat{x}_\theta)$ for different distributions of $c$.
These include the distribution used during training (i.e., $c \simiid {\cal U}([1, 2])$; leftmost boxplot), a distribution strictly contained within the training interval $[1,2]$, namely $[1.25, 1.75]$, and two intervals that progressively extend beyond the training range, namely $[0.75, 2.25]$ and $[0.5,2.5]$.
We note that the latter range implies that $\lambda \in [\pi,5\pi]$, so that there are still between 5 and 25 discretization points per wavelength.
The two figures reveal quasi identical behaviors and trends.
The best performance is observed when the network is used to solve systems involving matrices in the training domain, and the performance deteriorates as we move further away. 
However, for all problem sets, the numerical quality of the predicted solutions by the surrogate model is extremely poor, with a backward error $\eta_b(\hat{x}_\theta)$ that varies between $0.1$ and $0.8$.
Furthermore, the trend and performances are similar for the trained network; the results are not sensitive to the selected network.

\subsubsection{Generalization capabilities with respect to the domain size}\label{sec:raw_nn_gen_N}
\begin{figure}[!ht]\centering%
    \begin{subfigure}{.49\textwidth}%
    \includegraphics[width=\textwidth]{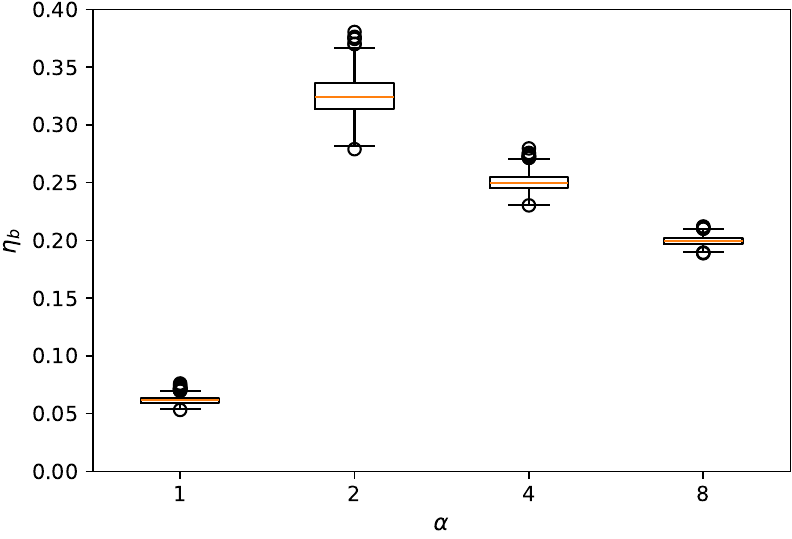}%
    \caption{Training \#1}%
    \end{subfigure}\hfill%
    \begin{subfigure}{.49\textwidth}%
    \includegraphics[width=\textwidth]{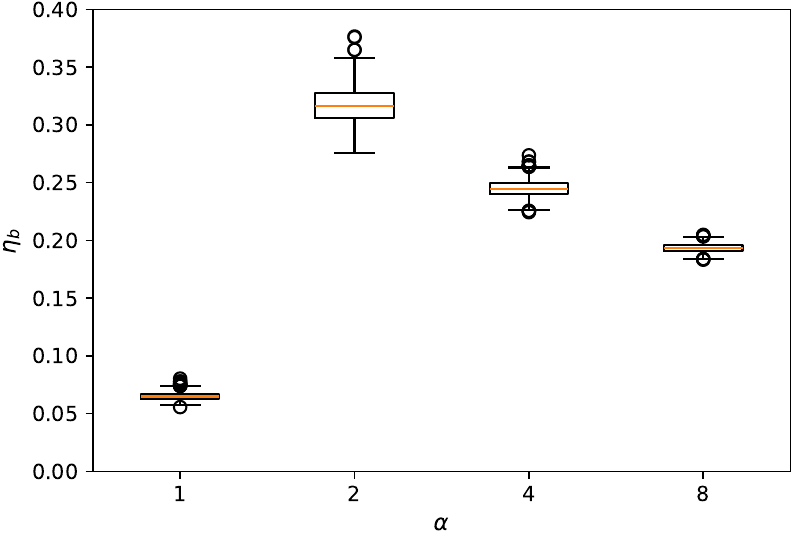}%
    \caption{Training \#2}%
    \end{subfigure}%
    \caption{Direct prediction from the neural network on problems with $\alpha \in \{1, 2, 4, 8\}$, where the neural network has been trained with $N=64$ and $c \simiid \mathcal{U}([1,2])$.
    Boxplots of $\eta_b$ based on \num{1000} experiments.}
    \label{fig:raw_nn_gen_dom_size}
\end{figure}

We conduct similar experiments with both trained networks while varying the domain size, keeping $c$ within the same interval as used for training.
As discussed in Section~\ref{sec:training}, the network is trained on a \( N \times N \) Cartesian grid, with $N=64$, discretizing the square domain $[-L,L]^2$, with $L=20$ and a PML size $\ell=4$.
We scale both $N$ and $L$ by a factor \( \alpha \in \{2, 4, 8\} \), so that $L/N=5/16$ remains fixed, while keeping the PML size $\ell=4$ and the distribution of $c \simiid \mathcal{U}([1,2])$ fixed.
This results in larger discretization matrices \( A(c) \in \C^{n \times n} \) with $n=N^2$, whose size thus grows quadratically with $\alpha$: $n=\num{16384}$ for \( \alpha = 2 \), $n=\num{65536}$ for \( \alpha = 4 \), and $n=\num{262144}$ for \( \alpha = 8 \); while $n = \num{4096}$ for the reference ($\alpha=1$).
As evidenced by Figure~\ref{fig:raw_nn_gen_dom_size}, the trends and performance remain similar.
The predicted solution exhibits extremely poor numerical quality.
At first glance, one might notice an improvement in quality as the domain size increases. However, the differences are not significant, and the overall quality remains very poor.
More notably, the variance appears to decrease with larger domain sizes, though we currently have no explanation for this behavior.

The results presented in the last two subsections indicate that the network cannot be used as a stand alone solver, as it lacks the robustness needed for accurate solutions. In the next section, we will explore its use as a preconditioner within a subspace linear solver (here, FGMRES), where an approximate inverse may be sufficient to accelerate convergence to an accurate solution.

\subsection{Numerical behavior of the trained network as nonlinear preconditioner}
\label{sec:nn-precond}

In this section, we report on the numerical effectiveness of the trained networks when used as a nonlinear preconditioner to speed-up the convergence of a subspace method, namely FGMRES. 
We first compare its performance with respect to an algebraic preconditioner in Section~\ref{sec:algebraic_comp}.
Then, in the following sections, we investigate the generalization capability with respect to the velocity field parameter $c$ and with respect to the domain size.

\subsubsection{Comparison with an algebraic preconditioner}\label{sec:algebraic_comp}
For the sake of comparison, we also report on the numerical performance of a purely algebraic implicit\footnote{implicit in contrast to explicit preconditioners based on approximate inverse techniques} preconditioner built as the (complete) LU factorization of an approximation $\sp{A}(c)$ of $A(c)$ that consists in keeping only the largest entries.
More precisely, for a given threshold $\epsilon>0$, the entries of $\sp{A}(c)$ are
\begin{equation}\label{eq:sparse A}
    \begin{cases}
        \sp{a}(c)_{ii} = a(c)_{ii},\\
        \sp{a}(c)_{ij} = a(c)_{ij} & \mbox{if } \vert a(c)_{ij}\vert > \epsilon \vert a(c)_{ii}\vert,\\
        0 & \mbox{otherwise}.
    \end{cases}
\end{equation}
We stress that this is different from an ILU preconditioner based on $A(c)$; in particular, the LU factors are not expected to be sparse.
We display in Table~\ref{tab:sparsified_A_1} the ratio of kept entries  in the preconditioner and relative distance between the matrix and the implicit preconditioner. We should mention that, due to the chosen discretization technique, the matrix $A(c)$ is quite populated with many nonzero entries  but usually does not need to be explicitly formed to solve the associated linear system by a subspace method. Among the large number of nonzero entries, many are of small magnitude which explain the very low nonzero ratio for a moderate distance between the two matrices.
As expected the lower the threshold $\epsilon$, the closer to the original matrix and the denser the preconditioner. Since the experiments are run on 100 examples, the quantities reported in Table~\ref{tab:sparsified_A_1} are the mean values, along with the associated standard deviations indicated between parentheses.

{\sisetup{exponent-mode=scientific, round-mode=figures, round-precision=2}%
\begin{table}[!ht]\centering%
    \begin{tabular}{S[table-format=2.2, exponent-mode=fixed, round-precision=3, round-pad=false]@{\qquad}c@{\qquad}c}%
      \toprule
      {$\epsilon$} & {$\nnz(\sp{A}(c)) / \nnz(A(c))$} & {$\|\sp{A}(c) - A(c)\|_\frob / \|A(c)\|_\frob$} \\
      \midrule
      0.1   & \num{3.411497e-03}~(\num{3.060292e-06}) & \num{1.803282e-01}~(\num{6.317130e-05}) \\
      0.05  & \num{6.480183e-03}~(\num{2.243749e-06}) & \num{1.169992e-01}~(\num{2.672865e-05}) \\
      0.02  & \num{1.458815e-02}~(\num{6.757878e-06}) & \num{7.725034e-02}~(\num{1.337933e-05}) \\
      0.01  & \num{3.222742e-02}~(\num{8.740236e-06}) & \num{3.825805e-02}~(\num{9.968938e-06}) \\
      0.005 & \num{4.892559e-02}~(\num{8.931253e-06}) & \num{1.790883e-02}~(\num{7.263261e-06}) \\
      \bottomrule
    \end{tabular}
  \caption{Sparsity of $\sp{A}$ in terms of proportion of non-zero~(denoted as nnz) entries kept compared to $A$ and relative distance of $\sp{A}$ to $A$ in terms of Frobenius norm, for different values of $\epsilon$.
    The average value over 100 experiments is reported, along with the standard deviation between parentheses.
    Note that  $\nnz(A(c)) = \num[round-mode=none, exponent-mode = fixed, fixed-exponent = 0]{7547926}$, and $\|A(c)\|_\frob = \num{1.110106e+03}~(\num{1.833749e-01})$.}%
  \label{tab:sparsified_A_1}
\end{table}%
}

\begin{figure}[!ht]%
    \centering%
    \begin{subfigure}{.5\textwidth}%
    \includegraphics[width=\textwidth]{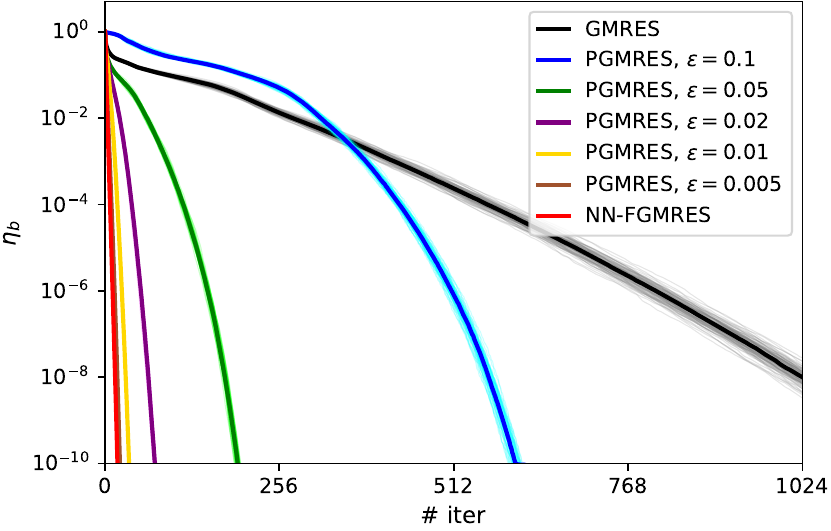}
    \caption{}
    \end{subfigure}\hfill%
    \begin{subfigure}{.49\textwidth}%
    \includegraphics[width=\textwidth]{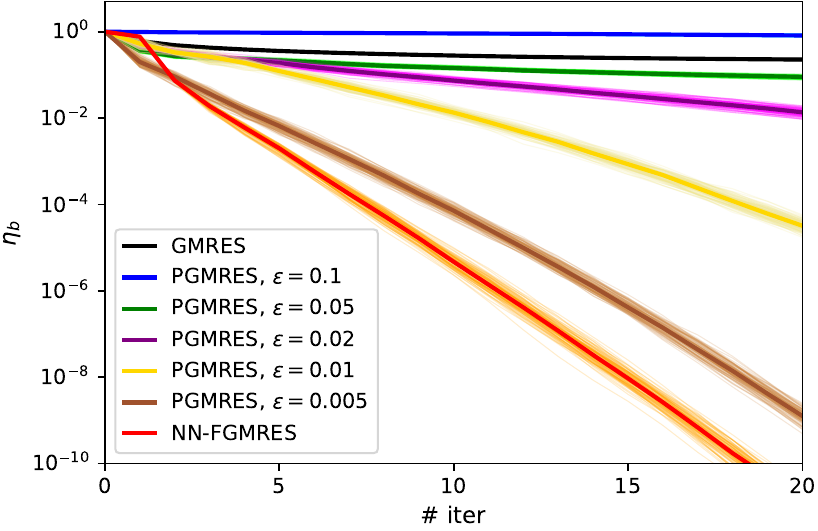}
    \caption{Zoom on the first 20 iterations.}%
    \label{fig:nn-fgmres-zoom}%
    \end{subfigure}%
    \caption{Convergence in terms of $\eta_b$ of unpreconditioned GMRES, preconditioned GMRES (PGMRES) with the algebraic preconditioners described in Table~\ref{tab:sparsified_A_1} and FGMRES preconditioned with the trained neural network (NN-FGMRES) corresponding to Training \#1, on problems with $\alpha=1$ and $c \simiid \mathcal{U}([1,2])$ (same as for the neural network training). 
    The iteration-wise median of the convergence history over 100 experiments is displayed as a thick line, while the individual 100 realizations are plotted as thinner, lighter lines.}
    \label{fig:nn-fgmres}
\end{figure}

Although the trained network performs poorly as a standalone solver, we highlight its potential to accelerate the convergence of subspace methods like FGMRES.
In Figure~\ref{fig:nn-fgmres}, we present the convergence history of $\eta_b$ for the solution of 100 linear systems derived from the discretization of the Helmholtz problem with a random velocity field $ c \simiid {\cal U}([1,2]) $, using FGMRES without restart. 
We display the convergence histories for the 100 linear systems, as well as the median for each system and for all preconditioners built using various sparsification thresholds $\epsilon$.
For the sake of completeness, we also display the convergence history of GMRES without a preconditioner, that  does not succeed to reach $10^{-10}$ in 1024 iterations while all the others do.
It can be observed that the algebraic implicit preconditioner becomes increasingly effective as $\epsilon$ decreases. 
It is noteworthy that $\NN$ consistently outperforms all the algebraic preconditioners, including the best-performing one obtained with $\epsilon = 0.005$, which retains nearly 5\% of the nonzero entries of each matrix and requires a new LU factorization for every linear system as $c$ varies.
In contrast, $\NN$ is the same for all linear systems---an important feature that makes $\NN$ particularly suitable as a preconditioner when one needs to solve a sequence of (possibly massively many) systems with different left- and/or right-hand sides. 
Interestingly, both approaches exhibit a similar spread around the mean performance, as highlighted in the zoomed-in view in Figure~\ref{fig:nn-fgmres-zoom}, suggesting that this variability is more intrinsic to the problem being solved than to the nature of the preconditioner itself.
Furthermore, all convergence histories using $\NN$ (in orange) are closely aligned and convergence to $10^{-10}$ is achieved in fewer than 20 iterations. These results illustrate the performance and the robustness of the trained network when used as preconditioner.

We conclude this section with a comparison with the complex shifted Laplacian (CSL) preconditioner~\cite{MagoluMongaMade2000_PreconditioningDiscreteHelmholtz,VanGijzen2007_SpectralAnalysisDiscrete}, which is widely used for Helmholtz-like problems.
CSL preconditioning aims to mitigate the mediocre behavior of ILU or multigrid (MG) preconditioners based directly on $A(c)$  due to the indefinite nature of the operator (see, e.g., \cite{Gander2005_IncompleteLuPreconditioner}).
The underlying principle of CSL preconditioning is to introduce a complex shift to the Laplace operator, in addition to the real-valued wavenumber shift of the Helmholtz equation~\eqref{eq:helmholtz}.
In the following experiment, we restrict ourselves to preconditioners of the form $M_{\beta}(c) \bydef A(c) + \jmath \beta k^2 I$, where $\beta \geq 0$, $k=\omega/c=1/c$, and $I$ denotes the identity matrix.
We note that $M_0(c) = A(c)$.
The tuning of $\beta$ is an active research question, and theoretical analyses are usually restricted to the case of uniform wavenumber fields $k$ (see, e.g., \cite{VanGijzen2007_SpectralAnalysisDiscrete,Gander2015_ApplyingGMRESHelmholtz}).
Generally speaking, small values of $\beta$ notoriously deteriorate the behavior of ILU- or MG-based CSL preconditioners (see, e.g., \cite{Gander2015_ApplyingGMRESHelmholtz}), so that practical values typically lie in the range $\beta \in [0.3, 1]$ (see, e.g., \cite{Erlangga2006_NovelMultigridBased}).
In what follows, we restrict ourselves to a similar range.

Because the operator $A(c)$ is dense in our Fourier-based discretization, it makes little sense to consider ILU factorizations.
Consequently, we consider again the \emph{complete} LU factorization of $M_{\beta}(c)$.
We note that this is not feasible in practice (otherwise, one could simply solve the discretized Helmholtz equation using a direct solver), but is merely used here for the sake of comparison.
In fact, the comparison is more than fair, as the complete LU decomposition of $M_{\beta}(c)$ is expected to behave better as a preconditioner than its ILU or MG counterpart, since it eliminates one source of approximation.

\begin{figure}[!ht]%
    \centering%
    \begin{subfigure}{.5\textwidth}%
    \includegraphics[width=\textwidth]{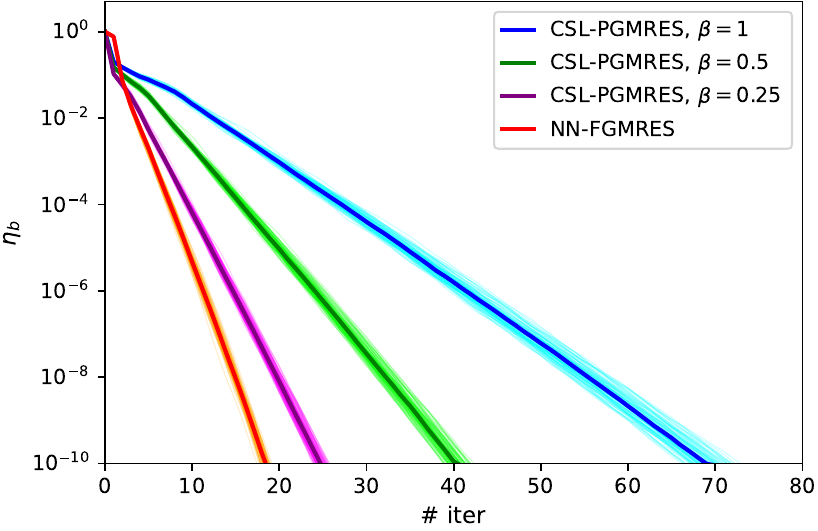}
    \caption{System-dependent CSL preconditioner.}
    \label{fig:nn-fgmres-csl-pb-dep}%
    \end{subfigure}\hfill%
    \begin{subfigure}{.49\textwidth}%
    \includegraphics[width=\textwidth]{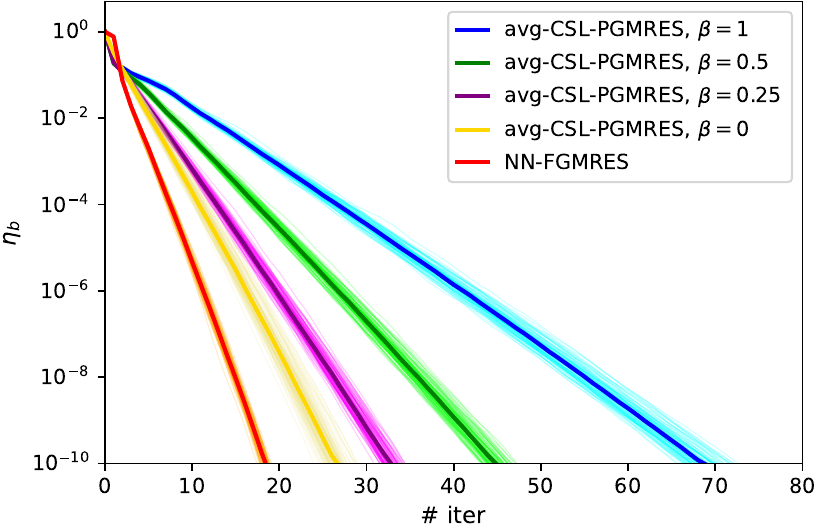}
    \caption{Constant CSL preconditioner based on $\bar{c}$.}%
    \label{fig:nn-fgmres-csl-avg}%
    \end{subfigure}%
    \caption{Convergence in terms of $\eta_b$ of the system-dependent (Fig.~\ref{fig:nn-fgmres-csl-pb-dep}, CSL-PGMRES) and average-based (Fig.~\ref{fig:nn-fgmres-csl-avg}, avg-CSL-PGMRE) CSL-preconditioned GMRES with various values of $\beta$, compared with NN-preconditioned FGMRES (NN-FGMRES) corresponding to Training \#1, on problems with $\alpha=1$ and $c \simiid \mathcal{U}([1,2])$ (same as for the neural network training). 
    The iteration-wise median of the convergence history over 100 experiments is displayed as a thick line, while the individual 100 realizations are plotted as thinner, lighter lines.}
    \label{fig:nn-fgmres-csl}
\end{figure}

Figure~\ref{fig:nn-fgmres-csl-pb-dep} shows the convergence history of the CSL-preconditioned GMRES, on the same problems as those of Figure~\ref{fig:nn-fgmres}, for various values of the parameter $\beta$ defining the purely imaginary shift of the CSL preconditioner $M_{\beta}(c)$.
As expected, the smaller $\beta$, the faster the convergence, which is expected as $M_{\beta}(c)$ gets closer to $A(c)$ (recall that $M_0=A$).
Nonetheless, for the smallest value of $\beta$ considered here, i.e., $\beta=0.25$, the CSL preconditioner is still outperformed by the NN preconditioner.
We stress again that the main advantage of the nonlinear NN preconditioner is that it is trained once and for all, and can subsequently be used to solve arbitrarily many problems.
On the contrary, the LU decomposition of $M_{\beta}(c)$ needed to be computed for every realization of $c$.
In practical applications, the same would apply for ILU factorizations; for MG-based CSL preconditioners, one or multiple MG cycles would need to be applied at every iteration of the outer subspace solver, for each individual problem defined by a realization of $c$.

A natural question is  whether constant CSL preconditioners may be designed.
To that end, we define a constant CSL preconditioner $M_{\beta}(\bar{c})$ based on the uniform, average discretized velocity field $\bar{c} \colon \vx \mapsto 1.5$.
We note that $M_0(\bar{c}) = A(\bar{c}) \ne A(c)$ for general realizations of $c$.
The convergence history of the constant CSL-preconditioned GMRES is shown is Figure~\ref{fig:nn-fgmres-csl-avg}.
For $\beta > 0$, the convergence is unsurprisingly slower than for the problem-dependent CSL preconditioner with the same value of $\beta$, owing to the additional approximation introduced by $\bar{c} \ne c$.
Even for $\beta = 0$, in this particular setting, the ``exact'' average-based preconditioner $M_0(\bar{c}) = A(\bar{c})$ is outperformed by the NN preconditioner.

In the next sections, we investigate the generalization capabilities of the NN preconditioner with respect to the parameter $c$ and with respect to the domain sizes.

\subsubsection{Generalization capability with respect to the parameter $c$}
As described in Section~\ref{sec:training}, the neural network is trained on a dataset built with $ c \in {\cal U}([1,2]) $ and demonstrates excellent performance as a preconditioner for solving linear systems whose matrices arise from the discretization of problems where $ c $ falls within the same range.  
A natural question arises: how does $ \NN$ perform on matrices constructed with $ c $ varying within an interval that is either entirely contained within or only partially overlaps the training interval $[1,2]$?  
To address this, we consider three sets of 100 linear systems each, with $ c $ varying in the following intervals: $[1.25,1.75]$, which is fully contained within the training range; $[0.75,2.25]$, where the training interval covers two thirds of the range; and finally $[0.5,2.5]$, where the training interval represents only one third of the values involved in the matrices preconditioned by $\NN$.

\begin{figure}[!ht]%
    \centering%
    \begin{subfigure}{.49\textwidth}%
    \includegraphics[width=\textwidth]{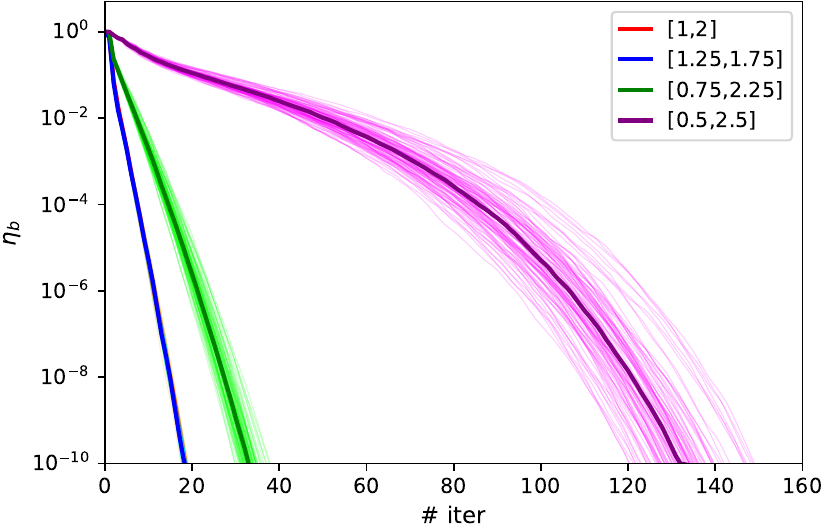}%
    \caption{Training \#1}%
    \end{subfigure}\hfill%
    \begin{subfigure}{.49\textwidth}%
    \includegraphics[width=\textwidth]{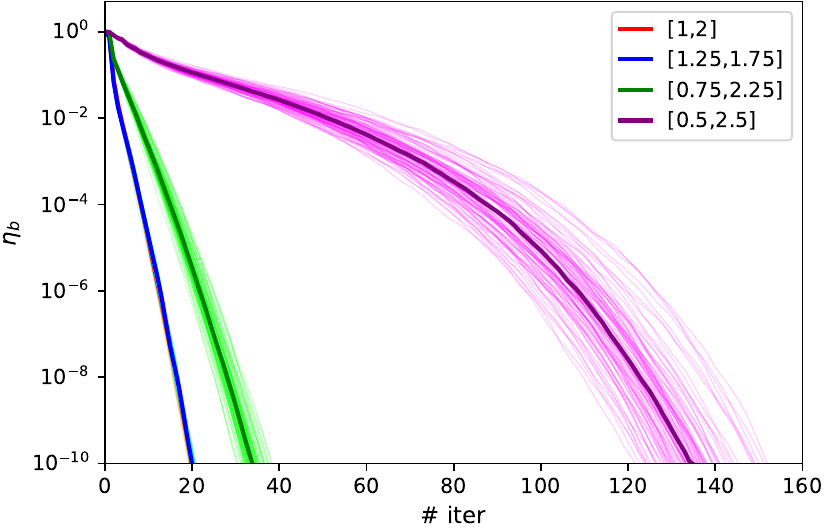}%
    \caption{Training \#2}%
    \end{subfigure}%
    \caption{Convergence in terms of $\eta_b$ of NN-FGMRES on problems with $\alpha = 1$ and $c \simiid \mathcal{U}([c_0,c_1])$ with varying ranges $[c_0,c_1]$ as reported, where the neural network has been trained with $\alpha=1$ and $c \simiid \mathcal{U}([1,2])$. 
    The iteration-wise median of the convergence history over 100 experiments is displayed as a thick line, while the individual 100 realizations are plotted as thinner, lighter lines.}
    \label{fig:nn-fgmres_gen_sos}
\end{figure}

The convergence history for these four sets of ranges (including the original range $[1,2]$) is shown in Figure~\ref{fig:nn-fgmres_gen_sos} for the two trained networks.  
The first key observation is that both trained networks exhibit very similar efficiency, indicating that the choice of the network does not significantly impact performance as long as they achieve comparable validation loss during training.  
The second observation is that when solving problems where $ c $ belongs to the training range or a subset of it, the convergence histories are nearly identical, as evidenced by the complete overlap of the curves (e.g., the red curves are entirely covered by the blue ones).  
Regarding generalization to values outside the training range, two insights emerge. First, the smaller the overlap with the training range, the greater the degradation in preconditioner performance---$\NN$ converges more slowly for $ c \in [0.5,2.5] $ than for $ c \in [0.75,2.25] $. Second, as the overlap decreases, the variability in the preconditioner’s effectiveness increases.  
Nonetheless, it is worth noting that despite the performance degradation, the preconditioned solver still achieves significant improvements compared to using no preconditioner at all.

We now further challenge the generalization capability of the nonlinear preconditioners by testing them on data (in terms of source term $b$ and velocity field $c$) with a structure that differs significantly from the training data.
Specifically, the entries of the discretized source term $b$ are set to zero except the entry corresponding to the center $(0,0)$ of the domain $\Omega$, which is set to one.
Furthermore, the velocity field $c$ is defined as a piecewise constant field with a constant value $c_0$ in the domain $\Omega$ except in a central, square domain $\square \bydef [-7.5, 7.5]^2 \subset \Omega$ where it has a constant value $c_1 > c_0$, in other words, $c(\cdot) = c_0 + (c_1-c_0) \indic_{\square}(\cdot)$, where $\indic$ denotes the indicator function.
The domain is kept to its original size, i.e., $\alpha=1$.
This setting will be referred to as the ``square'' setting.

\begin{figure}[!ht]\centering%
  \begin{subfigure}{.49\textwidth}\centering%
    \includegraphics[width=\textwidth]{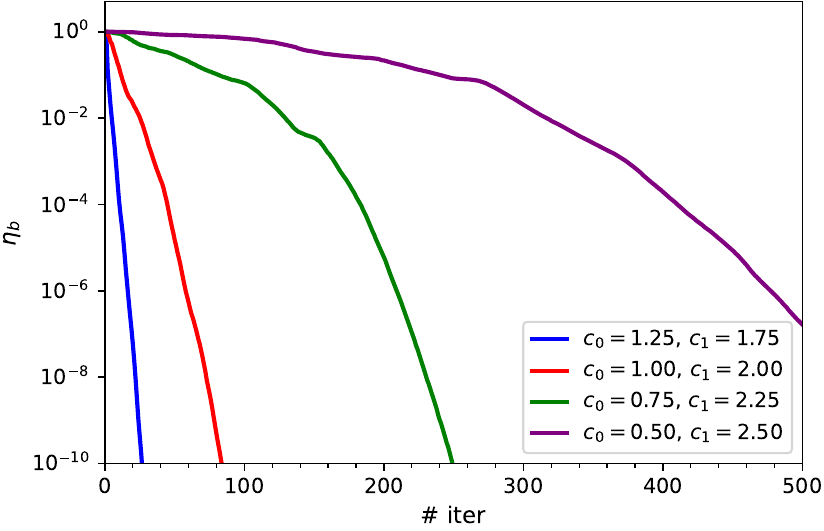}%
    \caption{Training \#1}
  \end{subfigure}\hfill%
  \begin{subfigure}{.49\textwidth}\centering%
    \includegraphics[width=\textwidth]{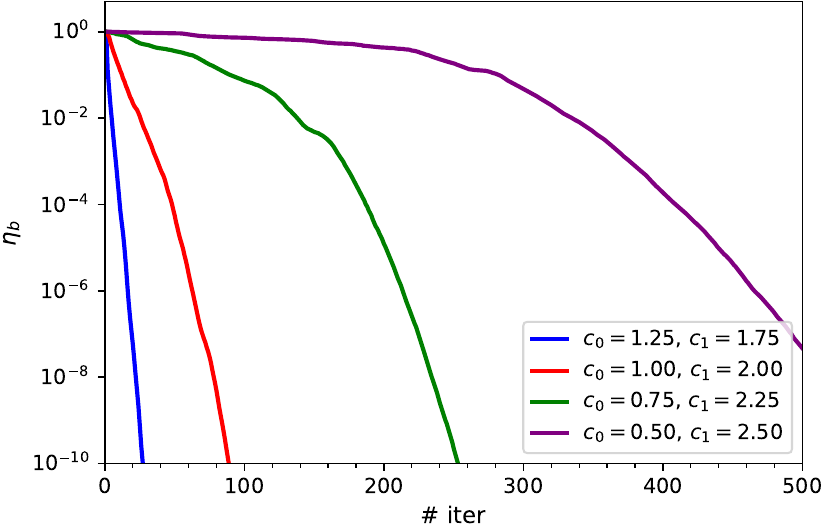}%
    \caption{Training \#2}
  \end{subfigure}%
  \caption{Convergence in terms of $\eta_b$ of NN-FGMRES on problems corresponding to the ``square'' setting, with different values of $c_0$ and $c_1$ as reported, where the neural network has been trained with $\alpha=1$ and $c \simiid \mathcal{U}([1,2])$.}
  \label{fig:nn-fgmres_gen_square}
\end{figure}

Figure~\ref{fig:nn-fgmres_gen_square} presents the convergence history of the NN-FGMRES solver for four ``square'' configurations, corresponding to four pairs $(c_0, c_1)$.
For the pair $(c_0 = 1.25, c_1 = 1.75)$, which fully lies in the training range $[1,2]$, the solver achieves fast convergence ($\eta_b$ reaches $10^{-10}$ in fewer than 30 iterations for both trainings), slightly slower than the reference (red curves in Figures~\ref{fig:nn-fgmres} and~\ref{fig:nn-fgmres_gen_sos}), despite the structures of the velocity field and source term being completely different from those seen during training.
For the pair $(c_0 = 1, c_1 = 2)$, which corresponds to the limit values of the training range, the convergence speed deteriorates but remains reasonable ($\eta_b$ reaches $10^{-10}$ in 85 to 90 iterations).
As $c_0$ and $c_1$ move further away from the training range, the convergence unsurprisingly keeps deteriorating.
Nevertheless, it seems that $\NN$ is still able to somewhat accelerate the convergence of the subspace solver.
These results suggest that the NN-FGMRES is highly robust to the PDE parameters $c$ and $b$, not only in terms of ranges of values, but also in terms of structure.

\subsubsection{Generalization capability with respect to the domain size}
\begin{figure}[!ht]\centering%
    \begin{subfigure}{.49\textwidth}\centering%
    \includegraphics[width=\textwidth]{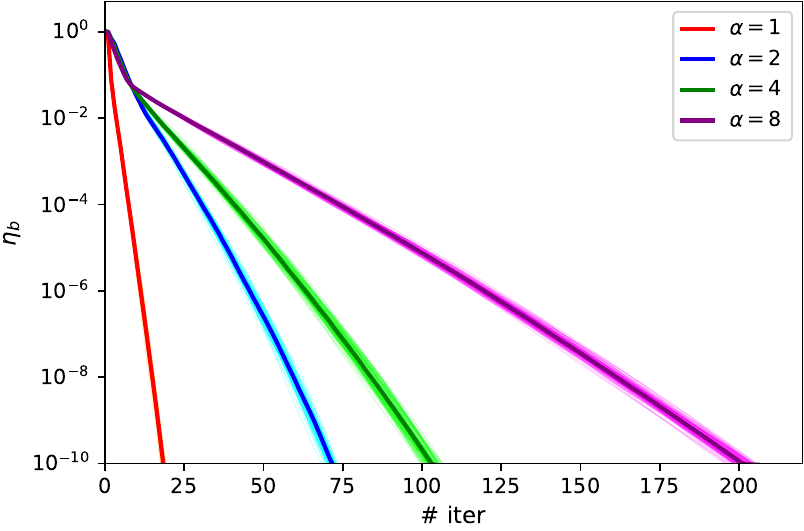}%
    \caption{Training \#1}
    \end{subfigure}\hfill%
    \begin{subfigure}{.49\textwidth}\centering%
    \includegraphics[width=\textwidth]{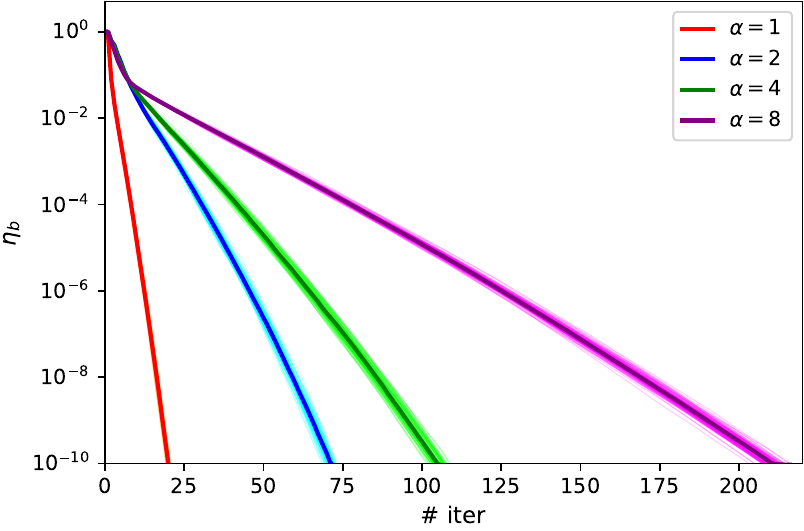}%
    \caption{Training \#2}
    \end{subfigure}%
    \caption{Convergence in terms of $\eta_b$ of NN-FGMRES on problems with $\alpha \in \{1, 2, 4, 8\}$, where the neural network has been trained with $\alpha=1$ and $c \simiid \mathcal{U}([1,2])$. 
    The iteration-wise median of the convergence history over 100 experiments is displayed as a thick line, while the individual 100 realizations are plotted as thinner, lighter lines.}
    \label{fig:nn-fgmres_gen_dom_size}
\end{figure}

In this section, we analyze the generalization capabilities of the preconditioner concerning the size of the physical domain, which translates into solving larger linear systems.  
Following the same approach as in previous experiments, we solve 100 linear systems for each system size using \( \NN \) with both trained networks. Figure~\ref{fig:nn-fgmres_gen_dom_size} displays the convergence history for all 100 solutions, along with the median convergence history.  
The results show that as the problem size increases, the number of iterations required for convergence also grows, though the dispersion around the median remains relatively stable. Additionally, the choice of trained network has little impact, as the trends for Training \#1 (left graph) and Training \#2 (right graph) are very similar.  
Notably, when increasing the scaling factor from \( \alpha = 4 \) to \( \alpha = 8 \), the number of iterations approximately doubles, even though the problem size is quadrupled. This highlights some robustness of the \( \NN \) preconditioner with respect to problem size.

Finally, we should mention that we deliberately chose not to report computational times, as our primary goal was to investigate the numerical behavior of the neural network-based preconditioner. No effort was made to optimize the implementation, and our study was limited to 2D cases. Nevertheless, due to the structure of the network, the inference phase at each iteration of the subspace method is extremely fast and would likely be difficult to outperform using classical scientific computing techniques related to any algebraic preconditioner.

\section{Further discussion}

We have presented results only with U-Net in this work, but preliminary results~\cite{Shpakovych2023_NeuralNetworkPreconditioning} on 1D examples suggest that using other neural network architectures,  trained with the same loss function leads to similar outcomes with different preconditioner efficiency. 
The objective of this work is to demonstrate that machine learning methods, which often struggle with precision issues, can be effectively combined with classical linear algebra techniques—such as suitable subspace methods—to leverage the best of both worlds. This approach ensures solution quality through traditional numerical techniques while benefiting from the acceleration provided by neural networks. It is a highly generic strategy that can be applied in various contexts~\cite{trifonov2025efficient}.
Identifying the optimal architecture is beyond the scope of this paper.  
However, certain machine learning architectures may lack the necessary properties to function effectively as nonlinear preconditioners. For example, as noted in~\cite{DeepONet}, DeepONet yields predictions that reside in a low-dimensional space whose dimension corresponds to embedding size (i.e., the output dimension of the trunk network). If this space has a lower dimension than that of the linear system to be solved, the resulting nonlinear preconditioner will be ``low-rank'', potentially leading to a breakdown of FGMRES. Therefore, careful consideration is required when designing the neural network architecture within the DeepONet framework to mitigate this issue. For an illustration of this potential limitation, we refer to~\cite{Shpakovych2023_NeuralNetworkPreconditioning}.

\section{Concluding remarks}
In this study, we illustrate  the potential of using neural networks as nonlinear preconditioners to solve parametric Helmholtz equations using subspace methods. While the accuracy and robustness of trained neural networks as standalone solvers is limited, integrating them as nonlinear preconditioners within the FGMRES algorithm significantly accelerates convergence, even for problems that are far beyond the training distribution. This hybrid approach combines the generalisation capabilities and computational speed of neural networks with the numerical reliability of classical iterative methods.

We would like to emphasise that the proposed strategy is not tied to a specific neural network architecture. Although we opted for a U-Net-based model due to its convolutional structure and efficient use of parameters, other architectures could also produce effective preconditioners, provided they retain sufficient representational capacity. The generic framework can be extended to other classes of linear systems beyond the Helmholtz equation, offering a promising approach to machine learning in scientific computing.

Several promising avenues for future research have emerged. Firstly, incorporating models that are independent of discretisation, such as PINNs or neural operators, which are trained across multiple resolutions, could improve the ability to generalise to arbitrary domain sizes and mesh structures. Secondly, extending the methodology to three-dimensional problems and other PDEs.

The intersection of numerical linear algebra and machine learning continues to pave the way for solving complex, high-dimensional scientific problems. The methodology presented here builds on this growing body of work by demonstrating that approximate learned inverses can be effectively integrated into established numerical schemes.

\section*{Acknowledgements}
We would like to thank Pierre Benjamin and Sofiane Haddad from Airbus CR\&T for fruitful discussions and exchanges that took place during the course of this work.
We are also grateful to Luciano Drozda and Selime Gürol from Cerfacs for their insightful feedback on an earlier draft of this work.

\printbibliography

\end{document}